\documentclass[11pt]{amsart2000}
\newtheorem{theorem}{Theorem}[section]
\newtheorem{lemma}[theorem]{Lemma}
\newtheorem{corollary}[theorem]{Corollary}
\newtheorem*{maintheorem}{Main Theorem}
\theoremstyle{definition}
\newtheorem{definition}[theorem]{Definition}
\numberwithin{equation}{section}

\begin{document}
\setlength{\unitlength}{0.01in}
\linethickness{0.01in}
\begin{center}
\begin{picture}(474,66)(0,0)
\multiput(0,66)(1,0){40}{\line(0,-1){24}}
\multiput(43,65)(1,-1){24}{\line(0,-1){40}}
\multiput(1,39)(1,-1){40}{\line(1,0){24}}
\multiput(70,2)(1,1){24}{\line(0,1){40}}
\multiput(72,0)(1,1){24}{\line(1,0){40}}
\multiput(97,66)(1,0){40}{\line(0,-1){40}}
\put(143,66){\makebox(0,0)[tl]{\footnotesize Proceedings of the Ninth Prague Topological Symposium}}
\put(143,50){\makebox(0,0)[tl]{\footnotesize Contributed papers from the symposium held in}}
\put(143,34){\makebox(0,0)[tl]{\footnotesize Prague, Czech Republic, August 19--25, 2001}}
\end{picture}
\end{center}
\vspace{0.25in}
\setcounter{page}{165}
\title{Compactification of a map which is mapped to itself}
\author{A. Iwanik}
\author{L. Janos}
\author{F. A. Smith}
\address{Department of Mathematics \& Computer Science\\
Kent State University\\
Kent, OH 44242}
\email{janos@mcs.kent.edu}
\email{fasmith@mcs.kent.edu}
\thanks{A. Iwanik, L. Janos and F. A. Smith,
{\em Compactification of a map which is mapped to itself},
Proceedings of the Ninth Prague Topological Symposium, (Prague, 2001),
pp.~165--169, Topology Atlas, Toronto, 2002}
\begin{abstract}
We prove that if $T:X\rightarrow X$ is a selfmap of a set $X$ such that
$\bigcap\left\{T^{n}X:n\in{\mathbb N}\right\}$ is a one-point set, then
the set $X$ can be endowed with a compact Hausdorff topology so that $T$
is continuous.
\end{abstract}
\maketitle

\section{Introduction}

If $(X,d)$ is a compact metric space and $T:X\rightarrow X$ is a Banach
contraction (there is $c\in[0,1)$ such that $d(Tx,Ty)\leq cd(x,y)$ for all
$x,y\in X)$ then the iterates $T^{n}$ shrink $X$ to a point $x^{\ast}$,
the unique fixed point of $T$. 
Thus $X$ and $T$ satisfy the set theoretical condition
\begin{equation}\label{eq1}
\bigcap\left\{T^{n}X:n\in{\mathbb N}\right\}=\left\{ x^{\ast}\right\}
\end{equation}

In the late 1960's J. deGroot asked about a converse to this, namely if
$X$ is an abstract set with cardinality at most that of the continuum and
$T$ a selfmap satisfying (\ref{eq1}) must there be a compact metric $d$
on $X$ so that $T$ is a Banach contraction? 
This and related questions have been examined see e.g.\ 
\cite{1,2,4,5,8}.
In \cite{4} the second author constructed a totally bounded metric
whenever (\ref{eq1}) is satisfied, but the question of compactness of $X$
remained open. 
In \cite{5} A.~Kub\v{e}na showed that a compact metric cannot exist by
constructing $2^{\mathfrak c}$ models of mutually nonisomorphic systems 
$(X_{i}, T_{i})$ satisfying (\ref{eq1}) and showing that the cardinality
of mutually nonisomorphic compactified systems $(X_{i}, \tau_{i}, T_{i})$
cannot exceed ${\mathfrak c}$.

As a result of this example, A. Iwanik, who only recently passed away,
asked if the conjecture of deGroot would be true with the metrizability
condition and the cardinality restriction removed. Indeed he,
together with the second and third authors, proved the following:

\begin{maintheorem}
If $X$ is a set and $T:X\rightarrow X$ satisfies the condition
(\ref{eq1}), then $X$ can be given a compact Hausdorff topology so that
$T$ is continuous.
\end{maintheorem}

In that which follows we provide a simplified proof of this theorem.

\section{Auxiliary Lemmas}

The main and dominating idea of compactification of sets respecting maps
between them comes to light by the following.

\begin{lemma}\label{lift}
If $X$ and $Y$ are disjoint sets and $T:X\rightarrow Y$ a surjective map,
then there are compact Hausdorff topologies on $X$ and $Y$ so that $T$ is
continuous.
\end{lemma}

\begin{proof}
Using the axiom of choice we well order the set $Y$ so that it has the
last element and endow $Y$ with the corresponding order topology. 
Thus $Y$ becomes a compact Hausdorff space. 
Doing the same with each set $T^{-1}y$ for $y\in Y$ we ``lift'' the
topology from $Y$ to $X$ ordering $X$ lexicographically according to the
order of $Y$.
Thus $X$ becomes also a compact Hausdorff space and since $T$ is evidently
order-preserving it is continuous.
\end{proof}

In the sequel we shall refer to a compact Hausdorff topology obtained by 
this method as a w.o.\ topology (well order topology).

This result and the technique in the proof will be applied systematically
many times in the sequel.
However, there is an obstacle to overcome. 
The maps we shall deal with are not surjective in general so that Lemma
\ref{lift} is not readily applicable. 
We must first partition the domains and targets of those maps into a
finite number of parts, called ``atoms'' so that the atoms will be mapped
onto atoms. 
We introduce some definitions concerning partitions of sets and their
behavior under mappings. 
By a partition $\pi$ of a set $X$ we mean a pairwise disjoint family of
sets $\{ C^{i}:i\in I\}$ such that 
$\bigcup\left\{C^{i}: i\in I\right\}=X$. 
If $\pi_{1}$, and $\pi_{2}$ are partitions then $\pi_{1}\leq \pi_{2}$
means that $\pi_{1}$ refines $\pi_{2}$ and $\pi_{1}\wedge \pi_{2}$ will
denote the common refinement of $\pi_{1}$ and $\pi_{2}$. 
If $T:X\rightarrow Y$ is a map between two disjoint sets $X$ and $Y$ and
if $\lambda$ is a partition of $Y$ defined by $\lambda=\{D^{j}:j\in J\}$
then $T^{-1}\lambda$ shall denote the partition of $X$ defined by
$\{T^{-1}D^{j}:j\in J\}$.

\begin{definition}\label{part} 
Let $T:X\rightarrow Y$ be a map between the disjoint sets $X$ and $Y$ and
let $\pi$ be a partition of $X$ given by $\pi=\{C^{i}: i\in I\}$. 
We denote by $T\pi$ the partition of $Y$ given by $T\pi=\{D^{j}: j\in J\}$
where $D^{j}$ are classes of the equivalence relation $\sim$ defined on
$Y$ by setting $y_{1}\sim y_{2}$ if the sets 
$\{i:T^{-1}y_{1}\bigcap C^{i}\neq \emptyset\}$
and $\{i:T^{-1}y_{2}\bigcap C^{i}\neq \emptyset\}$ coincide. 
It is evident that if $\pi$ is finite so is $T\pi$.
\end{definition}

\begin{definition}
If $X,Y$ and $T:X\rightarrow Y$ are as above and if $\pi$ and $\lambda$
are partitions of $X$ and $Y$, respectively, we say that $\pi$ and
$\lambda$ are $T$-{\em related} if every class of $\pi$ is mapped under
$T$ onto some class of $\lambda$. 
If $\pi$ and $\lambda$ are finite this is the desired situation. 
We say in this case that $\pi$ and $\lambda$ atomise $X$ and $Y$,
respectively.
\end{definition}

\begin{lemma}\label{related}
Let $T:X\rightarrow Y$ be as above and suppose that $\pi$ is a finite
partition of $X$ and $\lambda$ is a finite refinement of $T\pi$. 
Then the partitions $T^{-1}\lambda \wedge\pi$ and $\lambda$ are
$T$-related.
\end{lemma}

\begin{proof}
Let $\pi=\{C^{i}:i=1,\cdots n\}$, 
$T\pi=\{D^{j}:j=1,\cdots m\}$ and $\lambda=\{A^{k}:k=1,\cdots r\}$.

Every class of $T^{-1}\lambda\wedge\pi$ is a nonempty set of the form
$T^{-1}A^{k}\bigcap C^{i}$. 
This implies that there is some $a_{1}\in A^{k}$ with
\begin{equation}\label{eq2}
T^{-1}a_{1}\bigcap C^{i}\neq \emptyset.
\end{equation}
Since $\lambda$ is a refinement of $T\pi$ there is some $j$ for which
$A^{k}\subseteq D^{j}$. 
From this it follows that all elements of $A^{k}$ are equivalent under the
relation $\sim$ on $Y$ induced by $\pi$ (Definition \ref{part}). 
From this and \ref{eq2} it follows that 
$T^{-1}a\bigcap C^{i}\neq \emptyset$ for every $a\in A^{k}$ implying that
$T(T^{-1}A^{k}\bigcap C^{i})=A^{k}$ which concludes our proof.
\end{proof}

\begin{lemma}\label{compact}
Let $T:X\rightarrow Y$ be as above and suppose $\pi$ and $\lambda$ are 
finite $T$-related partitions of $X$ and $Y$, respectively. 
Then one can put on $X$ and $Y$ compact Hausdorff topology so that
\begin{itemize}
\item[(i)]
Each class of $\pi$ and $\lambda$ is compact
\item[(ii)]
$T$ is continuous.
\end{itemize}
\end{lemma}

\begin{proof}
Let $\pi=\{A^{i}:i=1,\cdots n\}, \lambda=\{B^{j}:j=1,\cdots m\}$.
We compactify $Y$ by putting on each set $B^{j}$ a w.o.\ topology and
for each $i=1, \cdots n$ we apply the techniques in the proof of 
Lemma \ref{lift} to lift this topology to a topology on $A^{i}$. 
Continuity of $T$ follows from continuity of each of its restriction to
the $A^{i}, i=1,\cdots n$.
\end{proof}

\begin{definition}
Let $X_{n}, n=1,\cdots, N, N\in{\mathbb N}$ be disjoint sets and
$T_{n}:X_{n}\rightarrow X_{n-1}, n=2, \cdots, N$ be maps. 
We shall call such finite family of sets and maps a chain of sets and
denote it by $\{X_{n}, T_{n}\}^{N}_{1}$.
\end{definition}

\begin{lemma}\label{chain}
For every chain $\{X_{n}, T_{n}\}^{N}_{1}$ there exist finite partitions
$\lambda_{1},$ $\ldots,$ $\lambda_{N}$ of $X_{1},$ $\ldots,$ $X_{N}$,
respectively
so that $\lambda_{n}$ and $\lambda_{n-1}$ are $T_{n}$-related for every
$n=2,\cdots, N$.
\end{lemma}

\begin{proof} 
We define inductively finite partitions $\pi_{n}$ of $X_{n}, n=1,\cdots N$
as follows. 
Starting with $\pi_{N}$ we set $\pi_{N}=\{X_{N}\}$ and if $\pi_{n}$ is
already defined we set $\pi_{n-1}=T_{n}\pi_{n}$. 
Thus $\pi_{1},\cdots \pi_{N}$ are finite partitions of 
$X_{1},\cdots, X_{N}$, respectively and we define the partition 
$\lambda_{1}$ of $X_{1}$ as $\pi_{1}$. 
We define $\lambda_{2}$ on $X_{2}$ as 
$T^{-1}_{2}\lambda_{1}\wedge \pi_{2}$ and if $\lambda_{n}$ is already
defined we define $\lambda_{n+1}$ as 
$T^{-1}_{n+1}\lambda_{n}\wedge \pi_{n+1}$. 
From the fact that $\lambda_{n}\leq \pi_{n}$ for $n=1,\cdots N$ and Lemma
\ref{related} we conclude that the partitions $\lambda_{n}, n=1,\cdots, N$
have the desired property, i.e., they ``atomize'' the chain
$\{X_{n},T_{n}\}^{N}_{1}$.
\end{proof}

From this result and Lemma \ref{compact} we obtain the following result:

\begin{theorem}\label{compactchain}
Any finite chain $\{X_{n}, T_{n}\}^{N}_{1}$ of sets and maps can be
compactified in the sense that one can put on each $X_{n}$ a compact
Hausdorff topology so that the maps $T_{n}:X_{n}\rightarrow X_{n-1},
n=2,\cdots N$ are continuous.
\end{theorem}

\begin{proof} 
By partitioning each $X_{n}$ by $\lambda_{n}$ as described by Lemma
2.4. we apply the argument used in the proof ofLemma \ref{chain}
sequentially to the maps $T_{2}, \cdots T_{n}$.
\end{proof}

\section{Proof of the Main Theorem}

Let $T:X\rightarrow X$ satisfy the condition (\ref{eq1}). 
This implies that the $T$-orbits $O(x)=\{T^{n}x:n\in{\mathbb N}\}$ are
either infinite or finite and in the latter case contain the fixed point
$x^{\ast}$. 
This allows us to visualize the system $(X,T)$ as a tree or more precisely
as a forest of trees. 
The individual trees will be defined as classes corresponding to the 
equivalence relation $\sim$ on the set $X\backslash \{x^{\ast}\}$ defined
by setting $x\sim y$ if there are $n\geq 0, m\geq 0$ such that
$T^{n}x=T^{m}y\neq x^{\ast}$. 
If $A$ is a class we say that $A$ is a class of the first kind if it
contains an element $z$ of the set
$T^{-1}x^{\ast}\backslash\{x^{\ast}\}$. 
In this case the class $A$ can be evidently represented as the disjoint
union of sets $\left\{ T^{-n}z: n\geq 0\right\}$; i.e.,
\begin{equation}\label{eq3}
A=\bigcup \left\{ T^{-n}z:n\geq 0\right\}
\end{equation}
We note that this family of sets $\{T^{-n}z:n\geq 0\}$ is finite since if
for every $n$ there existed a solution $x$ to the equation $T^{n}x=z$, 
(\ref{eq1}) would imply that $z=x^{\ast}$.

We now compactify the class $A$ by applying Theorem ref{compactchain}
to the chain $\{T^{-n}z, T_{n}\}$ where $T_{n}$ is defined as the
restriction of $T$ to $T^{-n}z$, $n\geq 1$.

Any class which is not of the first kind will be called one of the second
kind. 
If $A$ is of the second kind, let $a\in A$. 
Now $O(a)=\{b_{n}=T^{n}a:n\geq 0\}$ is infinite. 
For each $b_{n}\in O(a)$ let
\begin{equation}
B^{k}_{n}=\{x\in X:T^{k}x=b_{n}\ \mbox{and}\ T^{k-1}x\notin O(a)\}.
\end{equation}
For each fixed $n, \{B^{k}_{n}:k\geq 0\}$ forms a finite disjoint family,
since if not, (\ref{eq1}) implies $b_{n}=x^{\ast}$. 
We now apply Theorem \ref{compactchain} to the chain 
$\{B^{k}_{n}, T_{k}:k\geq 1\}$, where 
$T_{k}:B^{k}_{n}\rightarrow B^{k-1}_{n}$ is $T$ restricted to $B^{k}_{n}$. 
Thus, we compactify the set 
$B_{n}=\bigcup\left\{ B^{k}_{n}:k\geq 0\right\}$
which is the whole branch of the tree $A$ growing out of $T^{n}a$. 
The whole tree $A$ is now evidently the disjoint union of the sets $B_{n}$
which implies that $A$ receives a locally compact topology. 
Thus each class of the first kind is compact and each class of the second
kind is locally compact. 
This implies that the set $X\backslash\{x^{\ast}\}$ has a locally compact
topology as the disjoint union of locally compact spaces. 
We now compactify it by adding the point $x^{\ast}$. 
The map $T$ is continuous on $X\backslash \{x^{\ast}\}$ since its
restriction to each compact subset is. 
We now must show its continuity at $x^{\ast}$. 
This reduces to showing that for each open set $U$ containing $x^{\ast}$ 
the pre-image $T^{-1}U$ is again open. 
The set $U$ is of the form $X\backslash C$ where $C$ is a compact subset
of $X\backslash {x^{\ast}}$.
We observe by inspection easily that $T^{-1}C$ is again compact, so that
$T^{-1}(X\backslash C)=T^{-1}X\backslash T^{-1}C=X\backslash T^{-1}C$
which is again an open neighborhood of $x^{\ast}$ and which concludes the
proof of continuity of $T$, and with it the proof of our theorem.

\begin{corollary}
The deGroot conjecture is true if the set $X$ is countable.
\end{corollary}

\begin{proof}
This follows from the fact that a countable compact Hausdorff space is
metrizable.
\end{proof}


\begin{thebibliography}{1}

\bibitem{1}
J.~de~Groot and H.~de~Vries, \emph{Metrization of a set which is mapped into
  itself}, Quart. J. Math. Oxford Ser. (2) \textbf{9} (1958), 144--148. \MR{21
  \#4402}

\bibitem{8}
H.~de~Vries, \emph{Compactification of a set which is mapped onto itself},
  Bull. Acad. Polon. Sci. Cl. III. \textbf{5} (1957), 943--945, LXXIX.
  \MR{19,1069f}

\bibitem{2}
A.~Iwanik, \emph{How restrictive is topological dynamics?}, Comment. Math.
  Univ. Carolin. \textbf{38} (1997), no.~3, 563--569. \MR{98k:54078}

\bibitem{4}
Ludv{\'\i}k Jano{\v{s}}, \emph{An application of combinatorial techniques to a
  topological problem}, Bull. Austral. Math. Soc. \textbf{9} (1973), 439--443.
  \MR{49 \#3853}

\bibitem{5}
\bysame, \emph{Compactification and linearization of abstract dynamical
  systems}, Proceedings of the Eighth Prague Topological Symposium (1996)
  (North Bay, ON), Topology Atlas, 1997, pp.~157--162. \MR{99b:54067}

\end{thebibliography}
\providecommand{\bysame}{\leavevmode\hbox to3em{\hrulefill}\thinspace}
\providecommand{\MR}{\relax\ifhmode\unskip\space\fi MR }
\providecommand{\MRhref}[2]{%
  \href{http://www.ams.org/mathscinet-getitem?mr=#1}{#2}
}
\providecommand{\href}[2]{#2}

\end{document}